\documentclass{article}
\usepackage[dvips]{graphicx}
\usepackage{amsthm}  
\usepackage{amssymb}
\usepackage{amsmath}
\usepackage{enumerate, comment}
\usepackage{mathptmx}

\newcommand{\dt}{\circle*{3}}
\newcommand{\mrm}[1]{\mathrm{#1}}

\newcommand{\Aut}{\mrm{Aut}}
\newcommand{\ceil}[2]{\lceil \frac{#1}{#2}\rceil}

\theoremstyle{definition}
\newtheorem{defn}{Definition}[section]
\newtheorem{rmk}[defn]{Remark}

\newtheorem{note}[defn]{Note}

\newtheorem*{mainth}{Main Theorem}

\theoremstyle{plain}
\newtheorem{thm}[defn]{Theorem}

\newtheorem{lem}[defn]{Lemma}
\newtheorem{prop}[defn]{Proposition}

\begin{document}

\title{Maximally symmetric stable curves II}
\author{Michael A. van Opstall, R\u azvan Veliche}

\maketitle

\begin{abstract}
We determine a sharp bound for the automorphism groups of stable curves 
of genus $g$ with exactly $3g-3$ nodes.
\end{abstract}

\section{Introduction}

In a previous article \cite{msc}, we determined a sharp bound for the 
automorphism group
of a stable curve of genus $g$ and described all curves attaining the bound.
Aaron Bertram asked us if we could find the bound for those curves in the
zero dimensional strata of $\overline{M_g}$, that is, those stable curves
which have exactly $3g-3$ nodes (and thus have all rational components). We 
answer this question and also find
a sharp bound for the automorphism groups of such curves whose components are 
{\em smooth} and rational in this article. To state the result, we need a 
definition.

\begin{defn}
For a natural number $g$, define $b(g)$ to be the number of 1s in the binary
expansion of $g$, and $k(g)=g-b(g)$. Also set $h(g)=\lceil g/2\rceil -b(g)$.
\end{defn}

\begin{mainth}
Let $g\geq 2$ and let $C$ be a stable curve of genus $g$ with $3g-3$ nodes.
\begin{enumerate}
\item $|\Aut~C|\leq c\cdot 2^{g+k(g)}$. The value of $c$ for which this
bound is sharp is:
\begin{enumerate}
\item $c=3$ if $g=3\cdot 2^m$ for some $m\geq 0$.
\item $c=\frac32$ if $g=3\cdot 2^m+1$ for some $m>0$, or if $g=3(2^m+2^p)$ for
some $p\geq 0$, $m>p+1$.
\item $c=1$ otherwise.
\end{enumerate}
\item Suppose $g\geq 8$. If $C$ has only smooth components, then a 
bound for its automorphism group is given by $c\cdot 2^{g+h(g)}$. The value
of $c$ for which this bound is sharp is
\begin{enumerate}
\item $c=3$ if $g=3\cdot 2^m$ or $g=3(2^m+1)$ for some $m>1$
\item $c=\frac32$ if 
\begin{enumerate}
\item $g=3\cdot 2^m+1$ for some $m>1$,
\item $g=3\cdot 2^m+2$ for some $m>1$,
\item $g=3(2^m+2^p)$ with $p>0$ and $m>p+1$,
\item $g=3(2^m+2^p+1)$ with $p>0$ and $m>p+1$,
\end{enumerate}
\item $c=1$ otherwise.
\end{enumerate}
\end{enumerate}
\end{mainth}

The optimal graphs for $g<8$ in the second case do not follow a clear 
pattern, but can be found by exhaustive search. A graph whose automorphism
group has greater order than the bounds given above for $c=1$ will be called
{\em special}.

\section{Preliminaries}

Recall that a stable curve is a connected projective curve with at most nodes 
as singularities and finite automorphism group. The genus is computed as
$e-v+1$, where $e$ is the number of nodes and $v$ is the number of irreducible
components.

To each stable curve with rational components, one may associate a graph
knows as the {\em dual graph} of the curve.
The graph has one vertex for each component, and two vertices are joined by
an edge if and only if the corresponding components meet. If a component is
nodal, it meets itself, resulting in a loop in the graph. It is clear that
a stable curve of genus $g$ with $3g-3$ nodes has a trivalent dual graph
with $2g-2$ vertices.

Since the graphs are trivalent and the automorphism group of 
${\mathbb P}^1$is 
three-point transitive, the automorphisms of the curve are the same as the
automorphisms with the graph with a slight modification: we allow a
non-trivial automorphism to fix all the vertices in the graph. More concretely,
following what happens for the curves, a loop or a double edge contributes a
factor of two to the order of the automorphism group. The graph consisting of two
vertices joined by a triple edge has automorphism group of order twelve. 
Although this
is not the usual use of graph automorphism, it is the natural notion coming
from the problem, so it is what we will use. One may easily compute the 
order of the automorphism group of these graphs in the usual sense from our
formulas. From this point on, we will forget about the curve, and treat
the problem as a problem about graphs. Furthermore, all graphs will be 
assumed connected unless otherwise stated.

Our goal is to produce graphs which are as close to being trees as possible.
Define a {\em cone} as the graph which is a triangle with one edge doubled.
This graph is cubic except for one vertex of valence two, which will be used
to attach it to other graphs.

A tree has a unique vertex or edge which is common to all maximal length
geodesics. Call this the {\em root}. If $H$ is a graph with valence three
at all vertices but two, a {\em pseudocycle} of $H$s is a graph obtained
by replacing each vertex in a cycle by a copy of $H$.

Our proof will proceed by induction. We define the pinching of a simple
edge $e$ in a graph as follows: replace $e$ with two edges meeting in a
single vertex. This new vertex is valence two. We will use it to attach
the pinched graph to other graphs. Attaching will always be done in such
a way that the resulting graph is trivalent. For example, if two pinched
trivalent graphs are to be attached, a new edge must be introduced to 
join the pinched points.

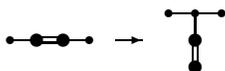
\begin{figure}[ht]
\begin{center}
\begin{picture}(90,30)
\put(5,15){\dt}
\put(5,15){\line(1,0){10}}
\put(15,15){\circle*{5}}
\put(15,14){\line(1,0){10}}
\put(15,16){\line(1,0){10}}
\put(25,15){\circle*{5}}
\put(25,15){\line(1,0){10}}
\put(35,15){\dt}
\put(45,15){\vector(1,0){10}}
\put(65,25){\dt}
\put(65,25){\line(1,0){10}}
\put(75,25){\dt}
\put(75,25){\line(1,0){10}}
\put(85,25){\dt}
\put(75,25){\line(0,-1){10}}
\put(75,15){\circle*{5}}
\put(74,15){\line(0,-1){10}}
\put(76,15){\line(0,-1){10}}
\put(75,5){\circle*{5}}
\end{picture}
\end{center}
\caption{Pinching double edges.}
\label{dep}
\end{figure}

Sometimes we will want to pinch a double edge. Here we will remove the double
edge and join the freed vertices. To this new vertex, attach a simple edge,
to which the double edge is attached. See Figure \ref{dep}. For completeness,
a pinched triple edge is a cone.

We refer to the star on four vertices simply as a ``star'', since we consider
no stars with more vertices. Loops and double edges will be considered 
cycles of length one and two in what follows. We will sometimes use the
word ``edge'' to refer to the graph which is two vertices joined by a single
edge. This should not be confusing. $O(e)$ will denote the orbit of an edge
$e$ in a graph $G$ under the action of the automorphism group of $G$.

Let $\Aut_e'G$ denote the group of automorphisms of the graph obtained by
pinching an edge (possibly multiple) $e$ of $G$. 

Define functions
\begin{eqnarray*}
\mu(g)&=&\max\frac{|\Aut~G|}{2^{g+h(g)}} \\
\mu_1(g)&=&\max\frac{|\Aut_e'G|}{2^{g+h(g)}},
\end{eqnarray*}
the maxima taken over all cubic graphs $G$ of genus $g$ and all edges $e$
in $G$. The value of $\mu(g)$ is the constant $c$ in the statement of the
Main Theorem. In fact, $\mu_1$ is always equal to one. Finally, define
$M(G)$ to be the number of edges in a minimal edge orbit in $G$.

At some stages in the proof, it will be most convenient to quote the
results of \cite{ts}. The method of proof used in that article is nearly
identical to this, so hopefully this lack of self-containment will actually
shorten the technical proofs and make this article easier to read.

\begin{defn}
For a natural number $g$, define the functions
\begin{eqnarray*}
l(g)&=&\min\{k:g=\sum_{i=1}^k a_i\cdot 2^{n_i}, a_i\in\{1,3\}\} \\
o(g)&=&g-l(g).
\end{eqnarray*}
\end{defn}

\begin{thm}\label{tsmain}
Let $G$ be a simple cubic graph with $g\geq 9$, where $g$ is
one more than the difference between the number of edges and number of
vertices.
\begin{itemize}
\item If $g=9\cdot 2^m+s$ ($s=0,1,2$) ($m\geq 0$) except $g=10,11,19,20,38$, 
then $|\Aut~G|\leq 3\cdot 2^{o(g)}$.
\item If $g=3\cdot 2^m+s$ ($s=0,1,2$) ($m\geq 2$), or $g=9(2^m+2^p)$ (with 
$|m-p|\geq 5$) or if $g=10,11,19,20,38$, 
then $|\Aut~G|\leq \frac{3}{2}\cdot 2^{o(g)}$.
\item If $g=5\cdot a\cdot 2^m+1$ (where $a=1$ or $3$, $m\geq 2$), then 
$|\Aut~G|\leq \frac{5}{4}\cdot 2^{o(g)}$.
\item Otherwise, $|\Aut~G|\leq 2^{o(g)}$.
\end{itemize}
Moreover, these bounds are sharp.
\end{thm}

\section{The candidates}

\begin{defn}
For every natural number $n>2$ we define a tree $T_n$ according to these
rules:
\begin{enumerate}
\item Place $n$ vertices in a row. This is called level one.
\item Assume level $k$ has been formed. In level $k+1$, place a vertex
and connect it to each of the first two vertices in level $k$. Continue
likewise until all pairs of level $k$ are exhausted. There may be a
vertex in level $k$ not connected to any vertex in level $k+1$.
\item If there is an unpaired vertex in level $k$ and one in some
level $l<k$, place a vertex in level $k+1$ and connect these two 
vertices to it.
\item If at some stage there are only two vertices left unpaired,
connect them with an edge.
\item If at some stage there are exactly three vertices left, connect all
three by edges to a new vertex.
\end{enumerate}
Define $T_2$ to be two vertices joined by an edge.
$T_n$ is trivalent at all interior vertices, and has $n$ leaves.
\end{defn}

\begin{figure}[ht]
\begin{center}
$\begin{array}{c@{\hspace{1in}}c}
\begin{picture}(50,30)
\put(5,25){\dt}
\put(15,25){\dt}
\put(25,25){\dt}
\put(35,25){\dt}
\put(45,25){\dt}
\put(5,25){\line(1,-2){5}}
\put(15,25){\line(-1,-2){5}}
\put(10,15){\dt}
\put(25,25){\line(1,-2){5}}
\put(35,25){\line(-1,-2){5}}
\put(30,15){\dt}
\put(10,15){\line(2,-1){20}}
\put(30,15){\line(0,-1){10}}
\put(30,5){\dt}
\put(45,25){\line(-3,-4){15}}
\end{picture}

&

\begin{picture}(70,30)
\put(5,25){\dt}
\put(15,25){\dt}
\put(25,25){\dt}
\put(35,25){\dt}
\put(5,25){\line(1,-2){5}}
\put(15,25){\line(-1,-2){5}}
\put(10,15){\dt}
\put(25,25){\line(1,-2){5}}
\put(35,25){\line(-1,-2){5}}
\put(30,15){\dt}
\put(10,15){\line(2,-1){20}}
\put(30,15){\line(0,-1){10}}
\put(30,5){\dt}
\put(30,5){\line(1,0){10}}
\put(40,5){\circle*{5}}
\put(40,4){\line(1,0){10}}
\put(40,6){\line(1,0){10}}
\put(50,5){\circle*{5}}
\put(50,5){\line(1,0){10}}
\put(60,5){\dt}
\end{picture}
\end{array}$
\end{center}
\caption{The tree $T_5$ and the core of $C_{11}'$.}
\label{t5}
\end{figure}
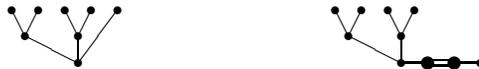

For example, this construction produces a binary tree when $n$ is a power
of two. Figure \ref{t5} shows the graph $T_5$.

The following proposition is a straightforward computation:

\begin{prop}\label{treeauts}
The order of the automorphism group of $T_n$ is $2^{k(n)}$, unless $n$
is of the form $n=3\cdot 2^m$, where the automorphism group has order
$3\cdot 2^{k(n)}$.
\end{prop}

\begin{defn}
For each $g>2$, define $C_g$ to be the graph constructed as follows:
\begin{itemize}
\item If $g\neq 3\cdot 2^m+1$ for $m>0$, and $g\neq 3(2^m+2^p)$ for $p\geq 0$,
$m>p+1$, $C_g$ is the graph formed by placing a loop at each tail of $T_g$.

\item If $g=3\cdot 2^m$, then the last stage in the construction of $T_n$
is connecting three vertices to a new vertex. $C_{g+1}$ is the graph formed by
following the construction of $T_n$ up until this last step, joining
these three vertices in a triangle rather than a star, and then placing a
loop on each tail.

\item If $g=3(2^m+2^p)$ for $p\geq 0$, $m>p+1$, attach a binary tree with 
$m$ leaves and a binary tree with $p$ leaves to the leaves of a binary
tree with two leaves. Arrange three copies of this configuration around a
star, and then place a loop on each tail. This graph is $C_g$.
\end{itemize}
\end{defn}

Again, it is easy to see:

\begin{prop}\label{graphauts}
The order of the automorphism group of $C_g$ is
\begin{enumerate}
\item $3\cdot 2^{g+k(g)}$ if $g=3\cdot 2^m$ for some $m\geq 0$,
\item $\frac32\cdot 2^{g+k(g)}$ if $g=3\cdot 2^m+1$ for some $m>0$ or if
$g=3(2^m+2^p)$ for $p\geq 0$, $m>p+1$,
\item $2^{g+k(g)}$ otherwise.
\end{enumerate}
\end{prop}

Now for the candidates for the case of all components smooth.

\begin{defn} Suppose $g>3$.
\begin{itemize}
\item If $g=3\cdot 2^m$ for some $m>0$, $C_g'$ is three binary trees with
$2^{m-1}$ leaves, connected with a star, with a cone attached to each leaf.
\item If $g=3\cdot 2^m+1$ for some $m>0$, $C_g'$ is obtained from 
$C_{g-1}'$ by expanding the central vertex to a triangle.
\item If $g=3\cdot 2^m+2$, $C_g'$ is obtained from $C_{g-2}'$ by expanding
the central vertex to a $K_{2,3}$.
\item If $g=3(2^m+1)$ for some $m>1$, $C_g'$ is obtained from
$C_{g-3}'$ by inserting a double edge in the middle of the edges of the
central star.
\item If $g=3(2^m+2^p)$ for $p>0$ and $m>p+1$ we construct $C_g'$ as follows:
attach two binary trees with $2^{m-1}$ and $2^{p-1}$ leaves at their roots.
Now attach three copies of this configuration to the ends of a star.
\item If $g=3(2^m+2^p+1)$ for $p>0$ and $m>p+1$, $C_g'$ is obtained from
$C_{g-3}'$ by inserting a double edge in the middle of the edges of the
central star.
\item Otherwise, 
\begin{itemize}
\item if $g$ is even, $C_g'$ is $C_{g/2}$ with the loops replaced
by cones;
\item if $g$ is odd, the formation of $T_{\lfloor g/2\rfloor}$ ends
with joining two vertices by an edge or joining three vertices by a star.
In the first case, insert a double edge in the center of the last edge
added, and place a cone at each leaf of the resulting graph to construct
$C_g'$. In the other case, since $g$ does not have any of the special
forms above, at least one of the graphs being attached to the star is not
isomorphic to the others. Insert a double edge in the center of the edge
of the star connected to this graph, place a cone at each leaf, and the
result is $C_g'$.
\end{itemize}
\end{itemize}
Define $C_2'$ to be a cone and $C_3'$ to be a cone with a double edge
attached to the free vertex by an extra edge. Note that $C_2'$ and $C_3'$ 
are not trivalent.

Although the graphs constructed in the ``otherwise'' part of this definition 
are not optimal in the exceptional genera, they have the important property 
that $\mu_1$ is one for these graphs.
Define $D_g$ for any $g$ to be the graph constructed in this last part.
\end{defn}

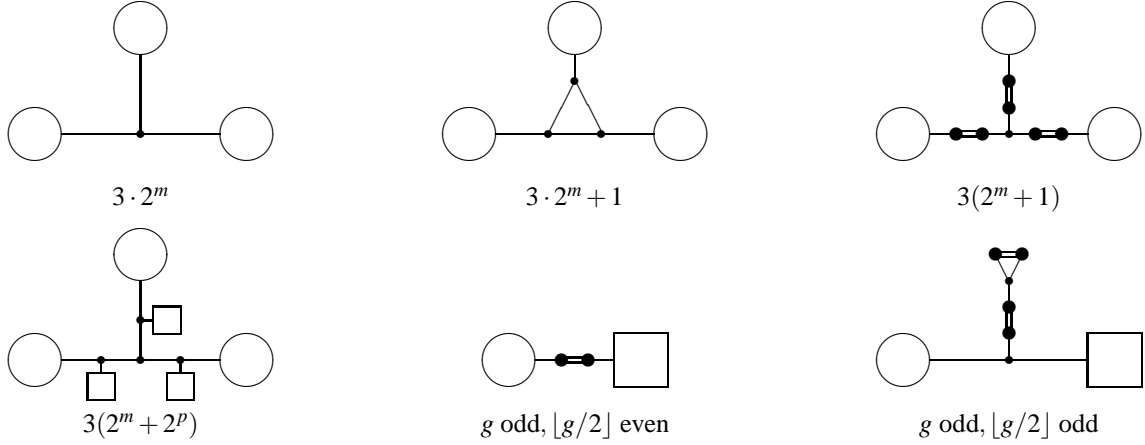
\begin{figure}[ht]
$\begin{array}{c@{\hspace{0.75in}}c@{\hspace{0.75in}}c}

\begin{picture}(110,70)
\put(15,15){\circle{20}}
\put(25,15){\line(1,0){60}}
\put(55,15){\dt}
\put(55,15){\line(0,1){30}}
\put(55,55){\circle{20}}
\put(95,15){\circle{20}}
\end{picture}
 &
\begin{picture}(110,70)
\put(15,15){\circle{20}}
\put(25,15){\line(1,0){60}}
\put(45,15){\dt}
\put(65,15){\dt}
\put(45,15){\line(1,2){10}}
\put(65,16){\line(-1,2){10}}
\put(55,35){\dt}
\put(55,35){\line(0,1){10}}
\put(55,55){\circle{20}}
\put(95,15){\circle{20}}
\end{picture}

&
\begin{picture}(110,70)
\put(15,15){\circle{20}}
\put(25,15){\line(1,0){10}}
\put(35,15){\circle*{5}}
\put(35,14){\line(1,0){10}}
\put(35,16){\line(1,0){10}}
\put(45,15){\circle*{5}}
\put(45,15){\line(1,0){10}}
\put(55,15){\dt}
\put(55,15){\line(1,0){10}}
\put(55,15){\line(0,1){10}}
\put(65,15){\circle*{5}}
\put(65,14){\line(1,0){10}}
\put(65,16){\line(1,0){10}}
\put(75,15){\circle*{5}}
\put(75,15){\line(1,0){10}}
\put(95,15){\circle{20}}
\put(55,25){\circle*{5}}
\put(54,25){\line(0,1){10}}
\put(56,25){\line(0,1){10}}
\put(55,35){\circle*{5}}
\put(55,35){\line(0,1){10}}
\put(55,55){\circle{20}}
\end{picture}

\\

3\cdot 2^m & 3\cdot 2^m+1 & 3(2^m+1) \\

\begin{picture}(110,70)
\put(15,15){\circle{20}}
\put(25,15){\line(1,0){60}}
\put(40,15){\dt}
\put(40,15){\line(0,-1){5}}
\put(35,0){\framebox(10,10){}}
\put(55,15){\dt}
\put(55,15){\line(0,1){30}}
\put(70,15){\dt}
\put(70,15){\line(0,-1){5}}
\put(65,0){\framebox(10,10){}}
\put(95,15){\circle{20}}
\put(55,30){\dt}
\put(55,30){\line(1,0){5}}
\put(60,25){\framebox(10,10){}}
\put(55,55){\circle{20}}
\end{picture}

&

\begin{picture}(80,30)
\put(15,15){\circle{20}}
\put(25,15){\line(1,0){10}}
\put(35,15){\circle*{5}}
\put(35,14){\line(1,0){10}}
\put(35,16){\line(1,0){10}}
\put(45,15){\circle*{5}}
\put(45,15){\line(1,0){10}}
\put(55,5){\framebox(20,20){}}
\end{picture}

&

\begin{picture}(110,70)
\put(15,15){\circle{20}}
\put(25,15){\line(1,0){60}}
\put(85,5){\framebox(20,20){}}
\put(55,15){\dt}
\put(55,15){\line(0,1){10}}
\put(55,25){\circle*{5}}
\put(54,25){\line(0,1){10}}
\put(56,25){\line(0,1){10}}
\put(55,35){\circle*{5}}
\put(55,35){\line(0,1){10}}
\put(55,45){\dt}
\put(55,45){\line(-1,2){5}}
\put(55,45){\line(1,2){5}}
\put(50,55){\circle*{5}}
\put(60,55){\circle*{5}}
\put(50,54){\line(1,0){10}}
\put(50,56){\line(1,0){10}}
\end{picture}

\\

3(2^m+2^p) & g\mrm{~odd}, \lfloor g/2\rfloor\mrm{~even} &
 g\mrm{~odd}, \lfloor g/2\rfloor\mrm{~odd}
\end{array}$
\caption{Schemata of the candidate graphs.}
\label{schem}
\end{figure}

Schematic diagrams for many of these cases are shown in Figure \ref{schem}.

The following two propositions are easy.

\begin{prop}
The order of the automorphism group of $C_g'$ is
\begin{enumerate}
\item $3\cdot 2^{g+h(g)}$ if $g=3\cdot 2^m$ or $g=3(2^m+1)$ for some $m>1$.
\item $\frac32\cdot 2^{g+h(g)}$ if $g=3\cdot 2^m+1$ or $g=3\cdot 2^m+2$ for 
some $m>1$.
\item $\frac32\cdot 2^{g+h(g)}$ if $g=3(2^m+2^p)$ or $g=3(2^m+2^p+1)$ with 
$p>0$ and $m>p+1$.
\item $2^{g+h(g)}$ otherwise.
\end{enumerate}
\end{prop}

\begin{prop}\label{growth.candidates}
$|\Aut~C_{g+1}|\geq|\Aut~C_g|$, with strict inequality unless $g=3\cdot 2^m$.
The same statement is true for $C_g'$.
\end{prop}

The proof of the Main Theorem for graphs with loops is similar to that for
graphs without loops, but easier. We will only prove the more difficult 
result here. However, the following proposition is useful in sorting out the
loop case, so we record it here.

\begin{prop}
If $G$ is a cubic graph of genus $g\geq 3$ with multiple edges, there exists a cubic graph
of the same genus without multiple edges whose automorphism group is at
least as large as that of $G$.
\end{prop}

\begin{proof}
There is a unique cubic graph with a triple edge. Its automorphism group has 
order twelve. On the other hand, there is a unique cubic graph on two vertices
without multiple edges. Its automorphism group has order eight. So we may
focus on double edges.

A double edge is adjacent to two simple edges. If these edges have a common
endpoint, we have a cone. Replace the cone with a binary tree with two leaves,
with a loop attached to each. This new configuration has strictly more 
automorphisms than the old.

If the two simple edges do not meet, then the double edge may be contracted
to a vertex, and a loop attached to an edge may be attached to this vertex, 
not affecting the automorphism group.
\end{proof}

\section{Small genus cases}

It follows from the results of our article \cite{ts} that the curves here have
more automorphisms than any curve of the same genus whose dual graph is
simple.

There is another series of graphs that is easily defined.
In genus $g$, define $C_g''$ as a $2(g-1)$-gon with every other edge double.
We have $|\Aut~C_g''|=2^g(g-1)$.

There are no simple cubic graphs on two vertices. A stable curve of genus 
two with all components smooth and rational is unique, corresponding to the
graph which is two vertices connected by a triple edge, for six automorphisms.

In genus three, the curve whose dual graph is a tetrahedron has 24 
automorphisms. The best that can be obtained with double edges corresponds
to the graph which looks like a square with one pair of opposite edges
doubled. This has only 16 automorphisms.

In genus four, the graph $K_{3,3}$ has 72 automorphisms. This is better
than $C_4'$ and $C_4''$, indeed, better than any non-simple graph.

In genus five, the optimal simple graph is a cube, which has 48 automorphisms.
This is beaten by $C_5'$ (64 automorphisms), which in turn is beaten by 
$C_5''$ (128 automorphisms).

In genus six, the Petersen graph is the optimal simple graph, with 120
automorphisms. It is beaten by $C_6''$ ($5\cdot 2^6$ automorphisms), which is
then beaten by $C_6'$ ($6\cdot 2^6$ automorphisms).

There is no question that simple graphs are always beaten after genus six.
In genus seven, $C_7''$ is the winner. In genus eight, $C_8'$ is optimal.
In genus nine, $C_9'$ and $C_9''$ have the same number of automorphisms.

The following table of values of these functions for low genus $g$
will be useful.

\begin{figure}[ht]
\begin{center}
$\begin{array}{c|c|c|c|c|c}
g & 2^{g+h(g)} & \mu & \mrm{graph~attaining~}\mu & \mu_1 & \mrm{graph~attaining~}\mu_1 \\ \hline
2 & 4 & 3 & \mrm{triple~edge} & 1 & \mrm{triple~edge} \\
3 & 8 & 3 & \mrm{tetrahedron} & 1 & C_3'' \\
4 & 32 & \frac94 & K_{3,3} & 1 & C_4' \\
5 & 64 & 2 & C_5'' & 1 & C_5' \\
6 & 128 & 3 & C_6' & 1 & D_6 \\
7 & 256 & 3 & C_7'' & 1 & D_7 \\
8 & 2048 & 1 & C_8' & 1 & C_8' \\
9 & 4096 & 1 & C_9', C_9'' & 1 & C_9'
\end{array}$
\end{center}
\caption{Table of small genus cases.}
\label{small.genus.table}
\end{figure}

\section{Structure of the proof}

For the remainder of the article, we will assume all results to be proved
hold for genus $h$ less than the genus $g$ under consideration. Call a double 
edge {\em free} if it does not occur as part of a cone.

\begin{lem}\label{structure.minimal.orbit}
Let $e$ be an edge of $G$ such that $O(e)$ has minimal order among all orbits 
of edges of $G$. Then only the following possibilities occur:
\begin{itemize}
\item $G=O(e)$;
\item $O(e)$ is a disjoint union of stars;
\item $O(e)$ is a disjoint union of edges;
\item $O(e)$ is a disjoint union of cycles; two such cycles are at distance 
at least two from each other.
\end{itemize}
\end{lem}

\begin{proof}
If is easy to see that if two stars in $O(e)$ have a common edge, then 
$G=O(e)$. Similarly, if two stars in $O(e)$ have a vertex in common, then 
either $G=O(e)$ or the third edge at that vertex will have an orbit of order 
smaller than that of $O(e)$ (which would be a contradiction to the choice of 
$e$).

Thus, if there is a star in $O(e)$, one of the first two possibilities occurs 
for $G$.

If no three edges in $O(e)$ share a common vertex, then either all edges in 
$O(e)$ are disjoint, or there are two edges $e_1$ (which may be assumed to be 
$e$, as $O(e)$ is acted upon transitively by $\Aut~G$) and $e_2$ in $O(e)$ 
with a common vertex $v$. Denote by $f$ the third edge of $G$ at $v$; $f$ is 
then not in $O(e)$. Denote by $w$ the other end of $e$.

If $v$ and $w$ are not in the same orbit of $\Aut~G$, then we see that 
$|O(e)|=2|O(v)|>|O(v)|\geq |O(f)|$, so we reach a contradiction to the 
choice of $e$. If however, $w\in O(v)$, then the existence of a cycle made 
of edges in $O(e)$ is immediate. Moreover, since $f\notin O(e)$, these cycles 
are disjoint.

Note that $|O(f)|\leq |O(e)|=|O(v)|$, with equality if and only if the ends 
of $f$ are not in the same orbit; in particular two cycles in $O(e)$ cannot 
be at distance one from each other (the edge between them, necessarily in the 
orbit of $f$, would have both endpoints in the same orbit).
\end{proof}

\begin{note}\label{e-f}
If the fourth situation above occurs, we will actually choose the edge $f$ 
and work with it in the arguments that follow; this is possible since $O(f)$ 
is also minimal, and may only be either a disjoint union of stars, or a 
disjoint union of edges. $f$ (or more precisely, its orbit) in this case will 
be called  {\em well-chosen}.
\end{note}

\begin{defn}
A graph of genus $g$ is {\em optimal} if its automorphism group has maximal
order among all automorphism groups of genus $g$ graphs. A graph $G$ is 
{\em strictly optimal} if in addition, the minimal order of the orbit of 
a well-chosen edge is smallest among those of all such optimal graphs, and
among such graphs with smallest minimal orbit, $G$ has the fewest double
edges.
\end{defn}

We will need to make some modifications to graphs in the course of the
proof. When we perform such a modification to an edge, vertex, or more
generally, any subgraph, we repeat the same construction at every other
subgraph in the order of the original. This ensures that we lose as few
automorphisms as possible.

In what follows, given a graph $G$, we will choose (well) an edge $e$ whose 
orbit
is minimal and remove the orbit, resulting in a graph $G'$. $G'$ will not
be cubic. There are two possibilities for the components of $G'$, since
$e$ is well-chosen: they are either of genus one or higher genus. For
higher genus components, there is a process of stabilization, which replaces
every path whose interior vertices have valence two and endpoints valence
three with a simple edge joining the endpoints. If the beginning and end
of the path coincide, this procedure could lead to loops. Since this takes
us outside of the class of graphs under consideration, it would be useful
to limit how often this can happen, so we can work inductively.

First, note that these loops happen only when $n-1$ endpoints of edges
in $O(e)$ are arranged around an $n$-cycle. It is not hard to see that $e$ 
could not have been minimal unless $n\leq 3$. The following proposition
limits the number of double edges that could lead to this problem.

\begin{lem}
A strictly optimal graph of genus $g\geq 10$ does not contain a path of 
alternating simple and double edges with two or more double edges.
\end{lem}

\begin{proof}
If there is a path of double edges alternating with
simple edges (some in a minimal orbit), let $k>1$ be the number of double 
edges in this path. If the beginning point and end point of this path are
the same vertex, either the whole graph is a pseudocycle, which is non-optimal
for $g\geq 10$, or the path begins and ends with simple edges. In this case,
remove the entire path and attach a copy of $C_k'$ to the free vertex 
resulting. If the path has distinct endpoints, replace
the path with a simple edge, then pinch this edge. To the pinched edge, 
attach a copy of $C_k'$. Usually, these two constructions increase the order 
of the
automorphism group. In any case, the attaching edge can be chosen as
an edge with minimal orbit. This new orbit is smaller than the original, 
unless there were exactly two double edges, in which case the new graph
has fewer double edges.
\end{proof}

\begin{prop}\label{fde}
There are fewer than eight free double edges in a strictly optimal 
graph $G$ of genus $g>9$.
\end{prop}

\begin{proof}
Suppose there are $k$ free double edges in $G$. Then form a new graph
$\bar{G}$ by replacing each configuration of a double edge with its two
adjacent edges by a single edge (this is lawful by the previous lemma).
$\bar{G}$ has genus $\bar{g}=g-k$ and the automorphism group has dropped
by a factor of $2^k$, so we have
\[
|\Aut~G|=2^k|\Aut~\bar{G}|.
\]
By induction (and the low genus table), the right hand side is smaller than 
$c\cdot 2^{g+h(\bar{g})}$ for some $c\leq 3$. Therefore, by our list of
candidates, $G$ was not optimal if $h(g)-h(\bar{g})\geq 2$. One may check 
that this is the case as soon as $k\geq 8$.
\end{proof}

The process of replacing double edges with simple ones (as in the above
proof) will be called {\em flattening}.

In cases of cycles of length three, one simply chooses the class of edges 
$f$ left adjacent to each triangle after removing $O(e)$. The orbit of $f$
is smaller or equal in size to that of $e$, and removing it and stabilizing
leads to double edges, but not loops.

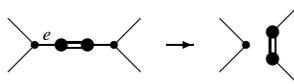
\begin{figure}[ht]
\begin{center}
\begin{picture}(120,40)
\put(5,10){\line(1,1){10}}
\put(5,30){\line(1,-1){10}}
\put(15,20){\dt}
\put(15,20){\line(1,0){10}}
\put(18,22){\makebox{\scriptsize $e$}}
\put(25,20){\circle*{5}}
\put(25,19){\line(1,0){10}}
\put(25,21){\line(1,0){10}}
\put(35,20){\circle*{5}}
\put(35,20){\line(1,0){10}}
\put(45,20){\dt}
\put(45,20){\line(1,1){10}}
\put(45,20){\line(1,-1){10}}
\put(65,20){\vector(1,0){10}}
\put(85,10){\line(1,1){10}}
\put(85,30){\line(1,-1){10}}
\put(95,20){\dt}
\put(105,15){\circle*{5}}
\put(104,15){\line(0,1){10}}
\put(106,15){\line(0,1){10}}
\put(105,25){\circle*{5}}
\put(105,15){\line(1,-1){10}}
\put(105,25){\line(1,1){10}}
\end{picture}
\end{center}
\caption{Stabilizing free double edges.}
\label{stabpic}
\end{figure}

\begin{rmk}\label{pinching}
When we remove a well-chosen minimal orbit of edges from a graph $G$ of the 
type we consider here, it is conceivable that when trying to stabilize the 
various components of $G'$ one is led to loops instead of just double edges. 
This happens in two cases only:
\begin{enumerate}
\item when a free double edge is incident at one vertex only to an edge in 
$O(e)$
\item when a free double edge is incident at both vertices to edges in 
$O(e)$; in this case at least an isomorphism class of components of $G'$ is 
made up of free double edges
\end{enumerate}
The second situation will be discussed separately, as it cannot be avoided 
even in the strictly optimal graphs.

However, we claim that the first situation is of no concern. Denote by $v$ 
and $w$ the vertices of a given free double edge, and assume that an edge 
$e$ in $O(e)$ is incident at $v$; denote by $u$ the vertex of $G$ at distance 
one from $w$, different from $v$; when $e$ is removed, enclose the double edge 
in the component it is part of in place of $u$ in the obvious way (see Figure
\ref{stabpic}) 
as to give a cubic graph. Except in a few cases (outlined below) the 
component such obtained is stable, and the estimates used in the ``regular'' 
cases (when no such manipulations were necessary) may be used to lead to the 
same restrictions on the number of edges in $O(e)$ and the structure of $G$ 
overall.

This construction should be done in a strictly optimal graph $G$; thus $u$ 
cannot be an end of a double edge by earlier remarks. However, it is possible 
that two or three double edges are adjancent to edges having $u$ as their
endpoint, and this could lead to conflicts. 

$u$ can be the merging point of three free double edges only when the component of $G'$ containing $u$ has genus three and is made up of a star with three double edges attached at its tails. Denoting by $s$ the number of these components, we see that removing them and joining the three edges in $O(e)$ to a common vertex produces a graph $\bar{G}$ with at most double edges (no triple edge is possible, as $e$ would be incident to two free double edges, which is not happening in a strictly optimal graph). Then by induction $|\Aut~G|\leq |\Aut~\bar{G}|\cdot 2^{3s}\leq 3\cdot 2^{g-3s+h(g-3s)+3s}\leq 2^{g+h(g)}$. We are interested in the last inequality being strict. Since $h(g)-h(g-3s)\geq h(3s)\geq 2$ for $s\geq 2$ (using \ref{inequalities}), we see that such a $G$ could not be optimal for $s\geq 2$. For $s=1$, we may change the minimal orbit to be that made up of the star itself, and then the above merging problems will not occur (since there are only three isolated edges incident to the new $O(e)$).

$u$ can be the merging point of two free double edges; in this case, the third edge at $u$ will have an orbit of order at most that of $O(e)$ (with equality possible only if $e$ sits between two free double edges, which should not happen in a strictly optimal graph $G$). Thus the strict optimality of $G$ and the minimality of $O(e)$ prohibit this situation.
\end{rmk}

\section{Proof of the Main Theorem}

We will repeatedly use an {\em exhausting subgraphs} argument. This entails 
choosing a connected component (star or edge) in $O(e)$, fixing its 
orientation (when the endpoints are in the same orbit) and then gradually 
enlarging the subgraph gotten at a certain stage by choosing one of its 
tails and adding whole components either of $O(e)$ or of $G'$ reached by that 
tail. When a component of $G'$ will be added, we will include in the new 
subgraph only the edges of $O(e)$ incident to it, and of these, in case 
$O(e)$ is a union of stars, only those that do not lead to stars whose 
center is already a vertex of the previous subgraph (in order to avoid 
cutting unnecessarily the number of tails).

At each step we look at the relative gain in the automorphism group. If a 
star is included at that step, then one of its edges is already fixed by the 
initial subgraph, then there could be at most a twofold increase in the order 
of the automorphism group at such a stage; moreover, such an increase occurs 
only when none of the vertices of the star was part of the subgraph at the 
beginning of the stage.

If however, a component is included at a certain step, then one of its 
vertices (which has valence two in $G'$) is already fixed, and that limits 
its symmetry; in other words, the automorphisms of the new subgraph fixing 
the previous one are precisely those fixing the incidence point. Once all 
these automorphisms are taken into account, all the edges incident to that 
component do not have extra freedom (they move where their incidence point 
moves), so may be added without further increase in the order of the 
automorphism group of the subgraph.

Unless otherwise noted, we will always expand the subgraphs by including 
whole components of $G'$ if the possibility exists (i.e. when not all tails 
of the subgraph gotten so far are centers of stars in $O(e)$). When the
exhausting finishes, we will have an estimate for the order of the automorphism
group of the graph, which can be compared to the order of the automorphism
group of the appropriate candidate graph.

We start by noting that:
\begin{itemize}
\item if $g=2u$, then $g+h(g)=3u-b(2u)=3u-b(u)$, while
\item if $g=2u+1$, then $g+h(g)=3u+2-b(2u+1)=3u+1-b(u)$
\end{itemize}

\begin{lem}\label{inequalities} The following inequalities hold:
\begin{enumerate}
\item $b(uv)\leq b(u)b(v)$; $b(u+v)\leq b(u)+b(v)$;
\item $h(u+v)\geq h(u)+h(v)$;
\item $h(uv+1)-uh(v)\geq \frac{u+1}{2}$ for $u\geq 4$;
\item $h(uv)-uh(v)\geq \lceil\frac{u-1}{2}\rceil$ for all $u\geq 1$, $v\geq 2$.
\end{enumerate}
\end{lem}

\begin{proof}
The proof is omitted, as the assertions are easy to check.
\end{proof}

\begin{lem}\label{stars.and.double.edges}
Assume that a graph $G$ has a minimal orbit of edges $O(e)$ made up of stars, 
and such that removing $O(e)$ leaves only cycles of length two (free double 
edges). Then $G$ cannot be optimal. 
\end{lem}

\begin{proof}
Let $k=|O(e)|$.

Flattening the free double edges and stabilizing the resulting graph leads 
to a graph $\bar{G}$; this is either of genus two, or is a simple cubic graph, 
which is moreover edge-transitive; moreover the genus $\bar{g}$ of $\bar{G}$ 
is $\frac{k}{2}+1$. This flattening implies 
$|\Aut~G|=2^k\cdot |\Aut~\bar{G}|$. From (\ref{tsmain}) we know that 
$|\Aut~\bar{G}|\leq 9\cdot 2^{\bar{g}-l(\bar{g})}$, so we get
$|\Aut~G|\leq 9\cdot 2^{\frac{3k}{2}+\frac{k}{2}+1-l(\frac{k}{2}+1)}=9\cdot 2^{2k+1-l(\frac{k}{2}+1)}$. We compare to $2^{g+h(g)}=2^{3k+1-b(k)}$. Then the inequality to study is $9\leq 2^{k+l(\frac{k}{2}+1)-b(k)}$, or, since $k=2u$ (even), $9\leq 2^{2u+l(u+1)-b(u)}$. Now $u\geq b(u)$ so the inequality is clearly strict for $u\geq 4$; one may check directly that the inequality is also strict for $u=2,3$. In case $u=1$, $k=2$ and $\bar{G}$ is a triple edge, and $g=5$; in this case $|\Aut~G|=12\cdot 2^3<8\cdot 2^4=|\Aut~C_5'|$, so again $G$ is not optimal.
\end{proof}

\begin{lem}\label{cycles.and.double.edges}
Assume that a graph $G$ of genus $g\geq 10$ has a minimal orbit of edges $O(e)$ made up of isolated (simple) edges, and such that $G'$ is made up of cycles only. Then $G$ is not optimal.
\end{lem}
\begin{proof}
Let $k=|O(e)|$.

If all cycles in $G'$ have length at least three, the graph is simple cubic, and (\ref{tsmain}) shows that this could not be optimal.

Consider first the case in which $G'$ has two isomorphism classes of cycles. 
 
If one class of components is made up of double edges (cycles of length two) and the other of cycles of length three and higher, then we may flatten the 
double edges and stabilize the resulting graph to $\bar{G}$; this is a simple cubic graph of genus $\bar{g}=g-k$; moreover, its minimal orbit of edges, 
has half the order of $O(e)$. We have $|\Aut~G|=2^k\cdot |\Aut~\bar{G}|$ and (\ref{tsmain}) shows that $\bar{G}$ could not have been optimal (for simple cubic graphs) anyway. Then we get $|\Aut~G|< 3\cdot 2^{g-l(g-k)}$ and would like to compare this to $2^{g+h(g)}$; we study then $3\leq 2^{\ceil{g}{2}+l(g-k)-b(g)}\geq 2^{1+\ceil{g}{2}-b(g)}$; this is easily seen to be strict for $g\geq 4$, and a moment's thought shows that such a graph cannot exist for $g\leq 3$.

The only case to consider is that in which all connected components (assumed to be cycles) of $G'$ are isomorphic of length two. Then $G$ must be a pseudocycle formed of $k$ double edges and $k$ simple edges joining them. Then $|\Aut~G|=2k\cdot 2^k$, $g=k+1$ and comparing to $2^{g+h(g)}$ leads to
$k\leq 2^{\ceil{k+1}{2}-b(k+1)}$, which is strict for $k\geq 9$ (with equality for $k=8$).
\end{proof}

\begin{rmk}\label{about.strictness}
The difference between strict and non-strict optimal graphs in genus $d$ becomes important only in genera $g>d$; more precisely, if two graphs have the same number of automorphisms and one has a minimal orbit with more edges than the other, their $\mu_1$ will be different (if need be we may use the lemma on the ratios of the orders of edges in a cubic graph to see that $\mu_1<1$ means $\mu_1\leq \frac{1}{2}$). In the induction step we determine which graph can reach or surpass the bound $2^{g+h(g)}$, assuming that $\mu_1=1$ for all components involved; as soon as one of these $\mu_1$ is strictly below one, the graph whose structure we analyze cannot possibly be optimal. Thus, whenever we will encounter graphs that could be optimal but not strictly optimal, we will check and mention explicitly the value of $\mu_1$ so that the next induction step is prepared.
\end{rmk}

\begin{prop}\label{stars}
If the minimal orbit of edges in $G$ is made up of stars, then it is made up 
of at most one star; moreover, the possible optimal graphs are only the 
$C_g'$ with $g=3\cdot 2^m$, $g=3(2^m+2^p)$, 
$g=3(2^m+2^p+1)$ and $g=3(2^m+1)$. In particular, the Main Theorem holds
in these cases. Furthermore, for these genera, $\mu_1(g)=1$.
\end{prop}

\begin{proof}
Let $k$ be the number of stars in $O(e)$.

We may stabilize the components of $G'$ according to
(\ref{pinching}), in light of the previous two lemmas. 

If $G'$ is connected, then 
$|\Aut~G|\leq |\Aut~G'|\leq 3\cdot2^{g'+h(g')}$; since $g'=g-2k\leq g-2$ 
and since $g+h(g)$ is strictly increasing (\ref{inequalities}) we get 
$3\cdot 2^{g'+h(g')}<2^{g+h(g)}$, so $G$ could not be optimal.

We concentrate then on the case when $G'$ is disconnected.

Let $t$ be the number of edges 
in $O(e)$ incident to a given component; we have $3k=st$. 
Using an expanding subgraphs argument, we get a bound for $|\Aut~G|$ as 
$6k\cdot 2^a\cdot 2^{s(d+h(d))}$  (assuming, by induction, that 
$\mu_1(d)\leq 1$ for all $d<g$) where $d$ is the genus of one of the $s$ 
isomorphic components of $G'=G\setminus O(e)$, and $a$ is the 
number of times we have no choice but to increase the subgraphs at a tail 
that is part of a star, and with the further restriction that the other edges 
in that star lead to components not incident to the subgraph obtained up to 
that point. We have $s\geq 3+2a$. Also, $g=s(d-1)+1+2k$. We are 
comparing the size of $|\Aut~G|$ to $2^{g+h(g)}$. Thus the inequality we are 
studying becomes:

$$6k\cdot 2^{\lfloor \frac{s-3}{2}\rfloor}\cdot 2^{s(d+h(d))}\leq 2^{g+h(g)}=2^{s(d-1)+1+2k+h(sd-s+2k+1)}$$
which is implied by
$$6k\leq 2^{2k+1-s+h(2k+1-s)+h(sd)-sh(d)-\lfloor\frac{s-3}2\rfloor}.$$

We remark that $t=1$ if and only if $k=1$ (a star and the components to which it is incident would give the whole connected $G$)
and that $s=1$ or $s\geq 3$ ($s=2$ would force a star to have at least two edges incident to a component, therefore all three edges would be incident to that component, therefore a component and the stars incident to it would give the whole graph; $G'$ would then be connnected, contradicting $s=2$).

We will simply sift the special cases through the filter of $2^{g+h(g)}$; each time a strict inequality is obtained, $G$ cannot be optimal, and we study the cases where the inequality fails separately.

If $t\geq 2$, then $s\leq \frac{3k}2$ and the inequality is implied by
$6k\leq 2^{\lceil \frac k2\rceil+h\left(\lceil\frac k2\rceil+1\right)+2}$
when $s\geq 4$; this is always a strict inequality. When $s=3$, easy 
manipulations of the desired inequality show that it is strict as long as 
$k\geq 3$. When $s=3$ and $k=2$ we have three connected components, each
incident to two stars. Working with the inequality shows that it fails
to hold only when $d=2^m$ $(m\geq 1)$. But then $g=3\cdot 2^m+2$, and 
$C_g'$ in these genera has a minimal orbit of three edges, so $G$ in this
case is not strictly optimal. In fact, such a $G$ is not even optimal, but
we do not need this, so the proof is omitted.

Therefore $t=1$, which implies $k=1$, so there is at most one star. In this
case, we have $s=3$, $v=1$, so the initial inequality becomes: 
$$6\leq 2^{3(d-1)+1+2+h(3d)-3d-3h(d)}=2^{d-2\lceil \frac{d}{2}\rceil +3b(d)-b(3d)}.$$

The last exponent is at least $b(d)-1$; so $b(d)\geq 4$ leads to strict inequality.

If $b(d)=3$ and $d$ is even, then again the inequality is strict; however, if $d$ is odd, then $b(d)=3$ and $b(3d)=6$ make the inequality fail. This may only happen when $d=2^m+2^p+1$ with $m>p+1>2$. In this case we obtain the exceptional $C_g'$ for $g=3\cdot (2^m+2^p+1)$.

If $b(d)=2$ and $d$ is even, the inequality is strict if $b(3d)\leq 3$, and fails for $b(3d)=4$. But $b(d)=2$ and $b(3d)=4$ for $d$ even force $d=2^m+2^p$ with $m>p+1>1$; this leads to the exceptional $C_g'$ for $g=3(2^m+2^p)$. If $b(d)=2$ and $d$ is odd, then the inequality is strict for $b(3d)\leq 2$, and fails for $b(3d)=3,4$. Then $d=2^m+1$ and $b(3d)=3$ cannot occur, while $b(3d)=4$ forces $m>1$; this leads to the $C_g'$ with $g=3(2^m+1)$.

If $b(d)=1$ and $d$ is odd, then $d=1$, which cannot happen for a cubic graph with a single star as a minimal orbit. If $d$ is even, then $d=2^m$ with $m\geq 1$ and then the inequality fails, revealing the graph $C_g'$ with $g=3\cdot 2^m$ for $m\geq 1$.

We note that all these graphs have $M=3$, and $\mu\leq 3$, therefore $\mu_1\leq 1$ (necessary for the induction step).
\end{proof}

\begin{prop}
If the minimal orbit of edges in a strictly optimal graph $G$ is made up of 
isolated (simple) edges, then 
\begin{enumerate}
\item There are at most three edges in that minimal orbit. 
\item $\mu_1(G)\leq 1$.
\item If $g=2^m$ or $g=2^m+1$, then $G$ must be $C_g'$.
of genus $g$ other than $G$.
\end{enumerate}
Finally, the Main Theorem holds in all cases not covered by the previous
proposition.
\end{prop}

\begin{proof}
By (\ref{pinching}), we may assume there are no problems stabilizing $G'$ when
$G'$ is connected.
Then $g'=g-k$ and  $|\Aut~G|\leq |\Aut~G'|\leq 3\cdot 2^{g'+h(g')}<2^{g+h(g)}$ 
for all $k\geq 2$, so $G$ could not be optimal. Moreover, when $k=1$, 
$\mu(G')\leq 1$ implies immediately that $G$ is not optimal, while 
$\mu(G')>1$ implies, by induction, that $M(G^\mrm{'stab})\geq 3$; but the 
marking inherited from stabilizing the two endpoints of the removed edge of 
$O(e)$ forces $M(G^\mrm{'stab})\leq 2$, so again $G$ could not be optimal.

We may now assume $G'$ is disconnected. 

{\bf Case I:} Suppose all the components of $G'$ have genus at least two.

{\bf Case I.1:} Further suppose $G'$ has two isomorphisms classes of 
components: one with 
$s_1$ components of genus $d_1\geq 2$, each with $n_1$ neighbours, each 
component connected to its neighbour by $t$ edges; similarly for 
$s_2,d_2,n_2$. We have $s_1n_1t=s_2n_2t=k$, $g=s_1d_1+s_2d_2+k+1-s_1-s_2$.

An expanding graph argument gives: 
$|\Aut~G|\leq k\cdot 2^{s_1(d_1+h(d_1))+s_2(d_2+h(d_s))}$; 
comparing this to $2^{g+h(g))}$ yields:

\[\tag{1}
k\leq 2^{k+1-s_1-s_2+h(s_1d_1+s_2d_2+k+1-s_1-s_2)-s_1h(d_1)-s_2h(d_2)}
\]

which is implied (using (\ref{inequalities})) by 

\[\tag{2}
k\leq 2^{k+1-s_1-s_2+h(k+1-s_1-s_2)+h(s_1d_1)-s_1h(d_1)+h(s_2d_2)-s_2h(d_2)}
\]

\begin{itemize}
\item If $t\geq 2$ and $\max(n_1,n_2)\geq 2$, then $\min(s_1,s_2)\leq \frac{k}{4}$ and $\max(s_1,s_2)\leq \frac{k}{2}$, also, $k\geq 4$; then the inequality (2) is implied by $k\leq 2^{1+\ceil{k}{4}+h(1+\ceil{k}{4})+1}$ (the last one in the exponent comes from the inequality $h(s_id_i)\geq s_ih(d_i)+1$ for at least one $i$). This is strict for all $k$.

\item If $t=1$ and $\min(n_1,n_2)\geq 3$ then $s_i\leq \frac{k}{3}$ so the inequality (2) is implied by $k\leq 2^{1+\ceil{k}{3}+h(1+\ceil{k}{3})}$ which is strict for all $k$.

\item If $t=1$ and $n_1=2<n_2$ (similarly for $n_2=2<n_1$) then $2\leq s_2\leq \frac{k}{3}, s_1\leq\frac{k}{2}$ and the inequality (2) is implied by $k\leq  2^{1+\ceil{k}{6}+h(1+\ceil{k}{6})+1}$ which is strict for all $k$.

\item If $t=1$ and $n_1=n_2=2$ then $2\leq s_i=\frac{k}{2}$ and writing $k=2u$ ($u=s_i\geq 2$) the inequality (2) is  implied by $2u\leq 2^{1+h(ud_1)-uh(d_1)+h(ud_2)-uh(d_2)}$. Since $d_1, d_2\geq 2$, this is implied by 
$u\leq 2^{\frac{u+1}{2}}$ which is strict for all $u$. 

\item If $t=1$ and $n_1=1$ then $s_1=k, s_2=1$ so the inequality (2) (or 
directly (1)) becomes 
$k\leq 2^{h(kd_1)-kh(d_1)}$ which is implied by
$k\leq 2^{\lceil\frac{k-1}2\rceil}$. The latter inequality is strict
for $k\geq 6$, while direct computation shows that the former is also strict
for $k=4,5$ and for all $k$ as soon as $b(d_1)\geq 3$. If $k=1$ the 
inequality cannot possibly fail.
\end{itemize}

Thus only the following subcases are left to analyze in detail:

\begin{itemize}
\item $t=2,s_1=s_2=1$ - one sees that this cannot occur for strictly optimal
$G$.
\item $t=1,s_2=1,s_1=k\leq 3, b(d_1)\leq 2$. When $k=2$, again, it is not
hard to see that $G$ cannot be strictly optimal. The case $k=3$ is actually
interesting. 
  \begin{itemize}
    \item If $b(d_1)\geq 3$, or $b(d_1)=2$ with $d_1$ even, strict 
inequality is obtained. 
    \item If $b(d_1)=2$ and $d_1$ is odd, $d_1=2^m+1$ 
$(m\geq 1)$. $G$ is made up of a central component of genus $d_2$, with three extremal components of genus $d_1=2^m+1$ attached each by an edge to the central component. Then we can bound $|\Aut~G|$ starting with the central component, so we get $c\cdot 2^{d_2+h(d_2)+3(d_1+h(d_1))}$ which compared to $2^{g+h(g)}$ yields $c\leq 2^{h(3d_1+d_2)-3h(d_1)-h(d_2)}$. The last exponent is $4+\ceil{d_2+1}{2}-\ceil{d_2}{2}+b(d_2)-b(3\cdot 2^m+3+d_2)$; as soon as it is at least two strict inequality is obtained. If $d_2=2u+1$ this exponent is $4+b(u)+1-b(3\cdot 2^m+4+2u)\geq 2$ so the inequality would be strict. Thus $d_2=2u$ and $u=1,3,4$ strict inequality follows, as it follows (for any $u$) when $c<2$. When $u=2$, $d_2=4$ and the only graph with $c\geq 2$ is the $K_{3,3}$; however, this is edge transitive and is pinched in three edges only, so it is easy to see that the corresponding $G$ is not optimal. In other cases, the Main Theorem gives inductively as the only possibility $d_2=3\cdot 2^p$, when $d_2\geq 10$; this immediately yields the structure of $G$ as the special $C_g'$. 
    \item If $b(d_1)=1$, since $d_1\geq 2$ we have $d_1=2^m$. Then argueing as in the previous case we get to the inequality $c\leq 2^{3+b(d_2)-b(3\cdot 2^m+d_2)}$; this is strict as soon as $c<2$. If $d_2\geq 10$ the Main Theorem yields as the only possibility for the central component $d_2=3\cdot 2^p$ or $3(2^p+1)$, which leads to the structure of $G$ as the $C_g'$ for the corresponding genus. In low genera, figure 4 shows that only $d_2=2$ may give an optimal structure (the pinching in three edges drops the number of automorphisms of the central component too much in the other cases). This yields the $C_g'$ for $g=3\cdot 2^m+2$.
  \end{itemize}
\end{itemize}

{\bf Case I.2:} Suppose there is a single isomorphism class of components of $G'$: there are $s$ components of genus $d\geq 2$, each with $n$ neighbours, linked to each of them by $t$ edges. $2k=snt$ and $g=sd+k+1-s$. 

As before, we get to 
\[\tag{3}
2k\leq 2^{k+1-s+h(sd+k+1-s)-sh(d)}
\]

implied by

\[\tag{4}
2k\leq 2^{k+1-s+h(k+1-s)+h(sd)-sh(d)}
\]

\begin{itemize}
\item $t\geq 2$, $n\geq 2$ implies $3\leq 3\leq \frac{k}{2}$ so the inequality (3) is implied by $k\leq 2^{\ceil{k}{2}+h(\ceil{h}{2}+1)+1}$ (the last one coming from $h(sd)-sh(d)\geq 1$ for $d\geq 2$ and $s\geq 2$); this is strict for all $k$.
\item $t\geq 2$, $n=1$ implies $s=2$ so $k=t$; the inequality (4) becomes $2k\leq 2^{k-1+h(k-1)+h(2d)-2h(d)}$; since $h(2d)-2h(d)\geq 1$ for $d\geq 2$, this inequality is strict except for $k=2$, and then only for $h(2d)-2h(d)=1$; this may only happen for $d=2^m$ or $d=2^m+1$ with $m\geq 1$. 
\item $t=1$, $n\geq 3$ implies $4\leq s\leq \frac{2k}{3}$ so the inequality (4) is implied by $2k\leq 2^{1+\ceil{k}{3}+h(1+\ceil{k}{3}))+1}$ which is strict for all $k$
\item $t=1$, $n=2$ implies $3\leq s=k$ ($G$ is a pseudocycle). The inequality (3) reduces to $k\leq 2^{h(kd+1)-kh(d)}$ which, as we have already seen, is strict except possibly for $b(d)=2$ and $k\leq 2$; since $k\geq 3$ we have actually strict inequality always
\item $t=1,n=1$ implies $s=2, k=1$ which reduces (3) to $2\leq 2^{h(2d+1)-2h(d)}$ which is strict except when $b(d)=2$ and $d$ is odd, or $b(d)=1$.
\item $t=2$, $n=1,s=2,k=2,d=2^m$ or $d=2^m+1$; these could be optimal, but easily seen not to be strictly optimal
\item $t=1,n=1,s=2,k=1$, $d=2^m$ or $d=2^m+1$; the latter does not lead to a strictly optimal $G$, while the first one gives the unique structure of the optimal $G$ in genus $2^{m+1}$.
\end{itemize}

This concludes consideration of $G'$ with all components of genus greater
than two.

{\bf Case II:} Now suppose $G'$ is disconnected and there are components 
which are cycles. Then (\ref{cycles.and.double.edges}) shows that 
we only need to deal with the case in which at least some components of $G'$ 
have genus two or more.

{\bf Case II.1:} If some components of $G'$ are free double edges (cycles of 
length two), 
then the edges of $O(e)$ are incident on both ends to these cycles; by strict 
optimality, the other components of $G'$ must stabilize properly. Then 
flattening the free double edges incident to $O(e)$ and stabilizing the 
resulting graph yields a new graph $\bar{G}$ which has at most double edges, 
genus $\bar{g}=g-k$ and with a minimal orbit of order half of $O(e)$; 
moreover, removing the minimal orbit of edges in $\bar{G}$ must leave only 
isomorphic components that stabilize properly (these the components of $G'$ 
which are not the cycles of length two). Then 
$|\Aut~G|\leq 2^k\cdot |\Aut~\bar{G}|$. Now 
$|\Aut~\bar{G}|\leq 2^{g-k+h(g-k)}$ implies $|\Aut~G|\leq 2^{g+h(g)}$, and if the first inequality is strict, the second is as well. Thus we only need to analyze the cases of equality or failure in the two Subcases 4 and 5 above and see what they mean for $G$. As before, let $s$ be the number of components which are not cycles, but have genus $d\geq 2$; two such components have $t$ edges 
incident to both in $\bar{G}$. Let $k'=\frac{k}{2}$; this is the order of the minimal orbit of $\bar{G}$.

The only cases to be considered are:

\begin{itemize}
\item $t=1, s=2, k=1, g-k=2^m$: this yields $G$ as the unique optimal graph for $g=2^m+1$ as two binary (cone tailed) trees with a double edge between their roots
\item $t=1,s=2,k=1, g-k=2^m+2$; this yields $G$ as two binary graphs of genus $2^m$ with their roots linked by a path containing three double edges. This is optimal, but not strictly optimal. As remarked before, this cannot appear as a component in constructing optimal graphs of higher genus, since its non-optimality would produce a $\mu_1\leq \frac{1}{2}$.
\end{itemize}

{\bf Case II.2:} If some components of $G'$ are cycles of length three or more, and the other components of genus $d\geq 2$ and stabilizing properly (i.e. the edges in $O(e)$ are not incident to free double edges) then denote by $s_1$ the number of cycles, by $s_2$ the number of other components; each cycle has $n_1$ neighbors among other components, and each other component has $n_2$ cycles as neighbors; denote by $t$ the number of edges in $O(e)$ incident to a fixed cycle and a fixed neighboring component. We have $k=s_1n_1t=s_2n_2t$, $g=s_2(d-1)+k+1$, and all cycle components of $G'$ have length $n_1t\geq 3$ (so $k\geq 3$).

Using an expanding subgraph argument, starting with a cycle and gradually adding up whole components of $G'$ as tails of the subgraph become incident to it, we get the estimate:

\[\tag{8}
|\Aut~G|\leq 2k\cdot 2^a\cdot 2^{s_2(d+h(d))}
\]

where $a\leq \lfloor \frac{s_2-n_1}{2}\rfloor$ represents the number of times a cycle is included in the expanding subgraph for lack of other choices (no other tail is incident to another component of $G'$).

Comparing this bound to $2^{g+h(g)}$ we get to:

\[\tag{9}
k\leq 2^{k-s_2-\lfloor \frac{s_2-n_1}{2}\rfloor + h(s_2d+k+1-s_2)-s_2h(d)}
\]

The last inequality in (\ref{inequalities}) now shows that this implied by 
$k\leq 2^{k-s_2-\lfloor \frac{s_2-n_1}{2}\rfloor+h(k+1-s_2)}$. 

\begin{itemize}
 \item If $t\geq 2$ and $n_2\geq 2$ (so $k\geq 4$), then $s_2\leq \frac{k}{4}$ and $a\leq \frac{k}{8}$ so the inequality (9) is implied by $k\leq 2^{\ceil{3k}{4}-\lfloor\frac{k}{8}\rfloor +h(\ceil{3k}{4}+1)}$ which is strict for $k\geq 4$
 \item If $t\geq 2$ and $n_2=1$, then $s_2\leq \frac{k}{2}, a\leq \frac{k}{4}$ so the inequality (9) is implied by $k\leq 2^{\ceil{k}{2}-\lfloor \frac{k}{4}\rfloor +h(\ceil{k}{2}+1)}$ which is strict for $k\geq 3$
 \item If $t=1$ and $n_2\geq 3$, then $s_2\leq \frac{k}{3}, a\leq \frac{k}{6}$ so the inequality (9) is implied by $k\leq 2^{\ceil{2k}{3}-\lfloor\frac{k}{6}\rfloor+h(\ceil{2k}{3}+1)}$ which is strict for $k\geq 3$
 \item If $t=1$ and $n_2=2$, then $s_2=\frac{k}{2}$, $k\geq 4$ even, $n_1\geq 3$ and $a\leq \frac{k-6}{4}$ so the inequality (9) is implied by $k\leq 2^{1+\frac{k}{2}-\lfloor \frac{k-2}{4}\rfloor +h(\frac{k}{2}+1)}$ which is strict for $k\geq 4$
 \item If $t=n_2=1$, then $k=s_2=n_1$ and the inequality (9) becomes $k\leq 2^{h(kd+1)-kh(d)}$. Then the second inequality in (\ref{inequalities}) shows that $k\geq 4$ and $b(d)\geq 2$ imply $h(kd+1)-kh(d)\geq h(kd)-kh(d)\geq \frac{k+1}{2}$; since $k<2^{\frac{k+1}{2}}$ for all $k\geq 1$ we get strict inequality in (9) as well. 
  \begin{itemize}
   \item If $b(d)\geq 3$ and $k=3$ then the inequality (9) becomes $3\leq 2^{\ceil{3d+1}{2}-3\ceil{d}{2}+3b(d)-b(3d+1)}\geq 2^{-1+b(d)-1}$; however, for equality to take place one needs $d$ to be odd and $b(3d+1)=2b(d)+1$, which is easily seen to be impossible
   \item If $b(d)=2$, $d=2p+1$ odd, $k=3$ then $h(3d+1)-3h(d)=-1+3b(d)-b(3d+1)\geq b(p)+1\geq 2$ so the inequality is strict.
   \item The case $b(d)=2$, $d$ even, and $k=3$ can be eliminated along the
lines of the arguments in subcase 4 of case 2.
   \item If $b(d)=1$ then $d\geq 2$ means $d=2^m,m\geq 1$. The inequality (9) becomes $k\leq 2^{k-b(k)}$ which is strict for $k\geq 4$, while for $k=3$ one obtains the $C_g'$ for $g=3\cdot 2^m+1$; again $\mu(G)=\frac{3}{2}$ and $M=3$.
  \end{itemize} 
\end{itemize}

We are done.
\end{proof}

\end{document}